\documentclass{article}[10pt]

\usepackage{
 amsmath,
 amsxtra,
 amsthm,
 amssymb,
 etex,
 mathrsfs,
  }
\usepackage[all]{xy}
\usepackage[T2A,T1]{fontenc}
\usepackage[OT2,T1]{fontenc}
\newcommand{\Zhe}{\mbox{\usefont{T2A}{\rmdefault}{m}{n}\CYRZH}}

\newtheorem{theorem}{Theorem}[section]
\newtheorem{lemma}[theorem]{Lemma}
\newtheorem{proposition}[theorem]{Proposition}
\newtheorem{corollary}[theorem]{Corollary}

\theoremstyle{definition}
\newtheorem{defn}[theorem]{Definition}
\newtheorem{conjecture}[theorem]{Conjecture}
\newtheorem{remark}[theorem]{Remark}

\newcommand{\bd}{\begin{defn}}
\newcommand{\ed}{\end{defn}}
\newcommand{\bl}{\begin{lemma}}
\newcommand{\el}{\end{lemma}}
\newcommand{\bp}{\begin{proposition}}
\newcommand{\ep}{\end{proposition}}
\newcommand{\bt}{\begin{theorem}}
\newcommand{\et}{\end{theorem}}
\newcommand{\bc}{\begin{corollary}}
\newcommand{\ec}{\end{corollary}}
\newcommand{\br}{\begin{remark}}
\newcommand{\er}{\end{remark}}
\newcommand{\ba}{\begin{array}}
\newcommand{\ea}{\end{array}}
\newcommand{\bpf}{\begin{proof}}
\newcommand{\epf}{\end{proof}}

\newcommand{\Z}{\mathbb{Z}}
\newcommand{\Q}{\mathbb{Q}}
\newcommand{\Zp}{\mathbb{Z}_{p}}
\newcommand{\Qp}{\mathbb{Q}_{p}}
\newcommand{\Op}{\mathcal{O}}
\newcommand{\Ep}{E[p^{\infty}]}
\newcommand{\Ap}{A[p^{\infty}]}

\newcommand{\Ga}{\Gamma}

\newcommand{\La}{\Lambda}
\newcommand{\la}{\lambda}

\DeclareMathOperator{\Sel}{Sel} \DeclareMathOperator{\Gal}{Gal}
\DeclareMathOperator{\Hom}{Hom} \DeclareMathOperator{\rank}{rank}
\DeclareMathOperator{\corank}{corank}
\DeclareMathOperator{\Ext}{Ext} 

\newcommand{\ord}{\mathrm{ord}}
\newcommand{\cyc}{\mathrm{cyc}}

\newcommand{\m}{\mathfrak{m}}
\newcommand{\mM}{\mathcal{M}}
\newcommand{\mG}{\mathcal{G}}
\newcommand{\mK}{\mathcal{K}}
\newcommand{\mL}{\mathcal{L}}
\newcommand{\mF}{\mathcal{F}}
\newcommand{\mT}{\mathcal{T}}
\newcommand{\mA}{\mathcal{A}}
\newcommand{\mR}{\mathcal{R}}

\newcommand{\ot}{\otimes}
\newcommand{\ilim}{\displaystyle \mathop{\varinjlim}\limits}

\newcommand{\coker}{\mathrm{coker}\,}

\newcommand{\lra}{\longrightarrow}

\newcommand{\ps}[1]{[[ #1 ]]}

  \DeclareFontFamily{U}{wncy}{}
  \DeclareFontShape{U}{wncy}{m}{n}{<->wncyr10}{}
  \DeclareSymbolFont{mcy}{U}{wncy}{m}{n}
  \DeclareMathSymbol{\sha}{\mathord}{mcy}{"58}

\begin{document}
\title{Structure of fine Selmer groups over $\Zp$-extensions}
 \author{Meng Fai Lim\footnote{School of Mathematics and Statistics $\&$ Hubei Key Laboratory of Mathematical Sciences, Central China Normal University, Wuhan, 430079, P.R.China. E-mail: \texttt{limmf@mail.ccnu.edu.cn}} }
\date{}
\maketitle

\begin{abstract} \footnotesize
\noindent This paper is concerned with the study of the fine Selmer group of an abelian variety over a $\Zp$-extension which is not necessarily cyclotomic. It has been conjectured that these fine Selmer groups are always torsion over $\Zp\ps{\Ga}$, where $\Ga$ is the Galois group of the $\Zp$-extension in question. In this paper, we shall provide several strong evidences towards this conjecture. Namely, we show that the conjectural torsionness is consistent with the pseudo-nullity conjecture of Coates-Sujatha. We also show that if the conjecture is known for the cyclotomic $\Zp$-extension, then it holds for almost all $\Zp$-extensions. We then carry out a similar study for the fine Selmer group of an elliptic modular form. When the modular forms are ordinary and come from a Hida family, we relate the torsionness of the fine Selmer groups of the specialization. This latter result allows us to show that the conjectural torsionness in certain cases is consistent with the growth number conjecture of Mazur. Finally, we end with some speculations on the torsionness of fine Selmer groups over an arbitrary $p$-adic Lie extension.


\end{abstract}

\section{Introduction}

The fine Selmer group has been ever-present in Iwasawa theory. Namely, it has been an object of frequent occurrence in the formulation (and proof) of the Iwasawa main conjecture (see \cite{K, Kob, PR00, Ru}). Despite this, it was only until the turn of the millennium that a systematic study of the said group was first undertook by Coates and Sujatha \cite{CS05a, CS05}, and a little later by Wuthrich \cite{WuPhD, Wu}, where they named it as we know today.
Initial studies mainly revolved around fine Selmer groups attached to abelian varieties (see loc. cit.; also see \cite{Hac, KL, LMu}). Subsequently, there have been much interest on the fine Selmer group of a modular form (for instance, see \cite{Jha, JS, HKLR}) or even more general classes of Galois representations (see \cite{LimPS, LimFine, LSu, KM}). A common feature in these cited works is that they are mainly concerned with working over the cyclotomic $\Zp$-extension.

Let $p$ be an odd prime. The goal of the paper is to consider the case of a $\Zp$-extension which is not the cyclotomic $\Zp$-extension. In this context, it has been conjectured that the fine Selmer group should be cotorsion over the ring $\Zp\ps{\Ga}$, where $\Ga$ denotes the Galois group of the $\Zp$-extension (see \cite{PR00, WuPhD, LimFineDoc}). This will be called Conjecture Y in the paper. One of our approaches towards studying this conjecture is guided by the following perspective which we now describe. Let $L_\infty$ be a $\Zp^2$-extension of a number field $F$, which is assumed to contain $F^\cyc$, the cyclotomic $\Zp$-extension of $F$. We shall write $G$ for the Galois group $\Gal(L_\infty/F)$. Let $A$ be an abelian variety over $F$. A conjecture of Coates-Sujatha \cite[Conjecture B]{CS05} (also see body of this paper) then predicts that the Pontryagin dual of the fine Selmer group, denoted by $Y(A/L_\infty)$, over the $\Zp^2$-extension $L_\infty$ is pseudo-null over the ring $\Zp\ps{G}$. Roughly speaking, this conjecture is saying that $Y(A/L_\infty)$ is ``quite small''. Therefore, in according to this conjecture of Coates-Sujatha, one would expect that the fine Selmer group at each $\Zp$-extensions of $F$ contained in $L_\infty$ should not be too ``large''. More concretely, we have the following observation.

\bp[Proposition \ref{torsion psuedo-null Zp2}] \label{torsion psuedo-null Zp2 intro}
 Let $A$ be an abelian variety defined over a number field $F$, and $L_{\infty}$ a $\Zp^2$-extension of $F$ which contains $F^\cyc$. Denote by $\Phi(L_\infty/F)$ the set of all $\Zp$-extensions of $F$ contained in $L_\infty$. Suppose that $Y(A/L_\infty)$ is pseudo-null over $\Zp\ps{\Gal(L_\infty/F)}$.
 Then $Y(A/\mL)$ is a torsion $\Zp\ps{\Gal(\mL/F)}$-module for every $\mL\in\Phi(L_\infty/F)$.
\ep

Although there are some numerical examples where pseudo-nullity of $Y(A/L_\infty)$ has been verified (for instance, see \cite{Jha, LeiP, LimFine}), the general statement remains wide open. The first of our main result is the following, which gives a way of obtaining infinite classes of $\Zp$-extensions, where Conjecture Y is valid.

\bt[Theorem \ref{torsion Zp2}] \label{torsion Zp2 intro}
 Let $A$ be an abelian variety defined over a number field $F$, and $L_{\infty}$ a $\Zp^2$-extension of $F$ which contains $F^\cyc$. Suppose that $Y(A/F^\cyc)$ is torsion over $\Zp\ps{\Gal(F^\cyc/F)}$. Then for all but finitely many $\mL\in\Phi(L_\infty/F)$, $Y(A/\mL)$ is torsion over $\Zp\ps{\Gal(\mL/F)}$.
 \et

The torsionness of $Y(A/F^\cyc)$ is known when $A$ is an elliptic curve over $\Q$ and $F$ is an abelian extension of $\Q$ (see \cite{K, Kob}). Hence the above theorem applies in these cases, where Conjecture Y is valid (also see Section \ref{examples} for some examples, where $F$ is not necessarily abelian over $\Q$ and the theorem applies). This therefore provides a strong evidence to Conjecture Y, and at the same time, a weak partial support towards the pseudo-nullity prediction of Coates-Sujatha.

As mentioned in the opening paragraph, there has been much interest in the study of the fine Selmer group of a modular form (for instance, see \cite{Jha, JS, HKLR}).
This will be the next theme of the paper which we describe briefly here.

Let $f$ be a normalized new cuspidal modular eigenform of even weight $k\geq 2$, level $N$ and nebentypus $\epsilon$. Write $A_f$ for the Galois module attached to $f$ (see body of the paper for its precise definition) which is defined over the ring of integers of $K_{f,\mathfrak{p}}$. Here $K_f$ is the the number field obtained by adjoining all the Fourier coefficients of $f$ to $\Q$, and $K_{f,\mathfrak{p}}$ is the localization of $K_f$ at some fixed prime $\mathfrak{p}$ of $K_f$ above $p$. Let $F_{\infty}$ be a $\Zp$-extension of $F$ with $F_n$ being the intermediate subfield satisfying $|F_n:F|=p^n$. We write $R(A_f/F_n)$ for the fine Selmer group defined over the field $F_n$ for $1\leq n\leq \infty$, and $Y(A_f/F_n)$ for its Pontryagin dual.
We then formulate an analogue of Conjecture Y for $Y(A_f/F_\infty)$, which will be called Conjecture Y$_f$ (see Conjecture \ref{Conj Yf}). Since an analogue of Conjecture B has been postulated and studied by Jha \cite{Jha}, we therefore have a natural analogue of Proposition \ref{torsion psuedo-null Zp2} for fine Selmer groups of modular forms (see Proposition \ref{torsion psuedo-null Zp2 modforms}). This in turn inspires the analogue of Theorem \ref{torsion Zp2} (see Theorem \ref{torsion Zp2 modforms}). The work of Kato \cite{K} again supplies many classes of examples, where $Y(A_f/F^\cyc)$ is torsion, thus allowing one to apply Theorem \ref{torsion Zp2 modforms} to obtain validity of Conjecture Y$_f$ in these cases.

We then follow up the above discussion by proving the following control theorem.

\bt[Theorem \ref{control theorem}] \label{control theorem intro}
Notations as above. Let $S$ be the set of primes of $F$ containing those dividing $pN$ and the infinite primes. Write $S_{fd}$ for the set of primes in $S$ which is finitely decomposed in $F_\infty/F$. Suppose further that either of the following statements is valid.
  \begin{enumerate}
    \item[$(i)$] For every prime $v\in S_{fd}$ and prime $v_n$ of $F_n$ above $v$,  $H^0(F_{n,v_n}, A_f)$ is finite.
    \item[$(ii)$] $K_{f,\mathfrak{p}}\cap \Qp(\mu_{p^\infty})=\Qp$ and $H^0(F_v, A_f)$ is finite for every prime $v\in S_{fd}$.
  \end{enumerate}
Then the restriction map
\[r_n: R(A_f/F_n) \lra R(A_f/F_{\infty})^{\Ga_n}\]
has finite kernel and cokernel which are bounded independently of $n$.
\et

For an abelian variety, a control theorem of such has been established by the author in \cite[Theorem 3.3]{LimFineDoc}. The above is therefore an analogue of this in the context of modular forms. Note that in this modular form context, there is an extra finiteness hypothesis on $H^0(F_v, A_f)$, and this arises due to a lack of an analogue of Mattuck's theorem \cite{Mat} for a modular form. We do however remark that although a recent work of Hatley-Kundu-Lei-Ray \cite{HKLR} has provided some sufficient conditions for this finiteness hypothesis to hold, it would seem that the general situation seems out of reach at the moment. We also note that in the event that the level $N$ is not divisible by $p$, then the finiteness is valid for all primes $v$ above $p$ (cf. \cite{CSW}).

We say a little more on the finiteness hypothesis on $H^0(F_v, A_f)$. As mentioned in the preceding paragraph, this is imposed on us by the lack of an analogue of Mattuck's theorem. In the proof of the control theorem, since we are estimating the kernel and cokernel at every intermediate $F_n$, the situation necessitates us to work al prior with a possible stronger hypothesis, namely, $H^0(F_{n,v_n}, A_f)$ is finite for every $v_n$ above $v$. As it turns out, in the event that $K_{f,\mathfrak{p}}\cap \Qp(\mu_{p^\infty})=\Qp$, the finiteness hypothesis at the base field suffices. In fact, we shall show that finiteness hypothesis at the base field $F$ will imply the finiteness hypothesis at every intermediate subfield $F_n$ (see Lemma \ref{rank leq 1 invariant} and proof of Theorem \ref{control theorem}). This latter observation seems interesting in its own right.

When the modular form arises from a specialization of an ordinary Hida deformation, we attach a fine Selmer group to the Hida deformation (denoted by $R(\mA/F_\infty)$; whose Pontryagin dual is denoted by $Y(\mA/F_\infty)$) and formulate an analogous conjecture which we call Conjecture $\mathcal{Y}$ (see Conjecture \ref{Conj mY}).  Our main result in this context is the following which can be thought as a ``horizontal'' variant of Theorems \ref{torsion Zp2} and \ref{torsion psuedo-null Zp2 modforms} (we refer readers to the body of the paper for the definitions of the objects and hypotheses appearing in the theorem).

\bt[Theorem \ref{Hida main}]
Assume that $\mathbf{(H1)}$ and $\mathbf{(H2)}$ are valid. Suppose that there exists $\eta\in\mathfrak{X}_{\mathrm{arith}}(h^{\ord}_{\mF})$ which satisfies the following properties.
\begin{enumerate}
  \item[$(a)$] For every prime $v\in S'_{cd}$, we have $H^0(F_{n,v_n}, A_f)$ being finite.
  \item[$(b)$] $R(A_{f_{\eta}}/F_\infty)$ is cotorsion over $\Op_\eta\ps{\Ga}$.
\end{enumerate}
Then Conjecture $\mathcal{Y}$ is valid, or equivalently, $R(\mA/F_\infty)$ is cotorsion over $\mR\ps{\Ga}$. Furthermore, for all but finitely many $\la\in\mathfrak{X}_{\mathrm{arith}}(h^{\ord}_{\mF})$, $R(A_{f_{\la}}/F_\infty)$ is cotorsion over $\Op_\la\ps{\Ga}$.
\et

Note that here again, the lack of an analogue of Mattuck's theorem necessitates us to work under assumption (a) (for a different set of primes).
We should mention that the above theorem is inspired by the work of Jha \cite{Jha}. Now, looking at Theorems \ref{torsion psuedo-null Zp2 modforms} and \ref{Hida main}, one can't help posing the following question.

\medskip \noindent
\textbf{Question $\mA$.} Is $R(\mA/F_\infty)^\vee$ pseudo-null over $\mR\ps{\Ga}$?
\medskip

Unfortunately, we do not have an answer to this. In fact, to the best knowledge of the author, even in the context of a cyclotomic $\Zp$-extension $F^\cyc$, the structure of $R(\mA/F^\cyc)$ does not seem well-understood (but see a very recent work of Lei-Palvannan \cite{LeiP2} for some discussion in this direction).

It should be evident to the readers that much of the discussion in this paper may be extended to fine Selmer groups attached to even broader classes of Galois representations of interest as studied in \cite{LimFine, LSu, KM}. We have decided to restrict our attention to the context of the paper to simplify the presentation. Furthermore, we believe that even in the modular form or Hida deformation context, the occurrence of certain interesting phenomenon deserves further future studies. (For instance, the lack of an analogue of Mattuck's theorem definitely requires further investigation and this sort of issues will naturally come up if one wants to study fine Selmer groups of more general  Galois representations.)

Although the focus of our paper is to formulate variants of Conjecture Y, we should remark that it would be interesting to study the variation of the Iwasawa invariants of the $\Zp$-specializations (either horizontally or vertically). We will not pursue this here but refer readers to \cite{G73, Cu, Mon, Kl} for some discussion in the vertical direction.
While we have nothing to say about this variational aspect, we shall end by formulating and investigating a generalized Conjecture Y (see Conjecture \ref{General Conj Y}) over an arbitrary $p$-adic Lie extension of dimension $>1$, and give some conceptual evidences towards the paucity of this generalized conjecture (see Remark \ref{K-groups remark} and Proposition \ref{torsion psuedo-null non-com}).

We now give an outline of our paper. In Section \ref{com alg review}, we collect several results on modules over regular local rings which will be required in our arithmetic discussion. In Section \ref{fineAbvar}, we introduce the fine Selmer groups of abelian varieties and establish Theorem \ref{torsion Zp2}. Section \ref{finemodform} is where we study the fine Selmer group of a modular form. The control theorem for the fine Selmer group of a modular form will be proved in this section. In Section \ref{ordinary sec}, we investigate the relationship between the Conjecture $\mathcal{Y}$ on the fine Selmer group
for a Hida family and the corresponding Conjecture Y$_f$ for the specializations. We also discuss a situation showing that the conjectural torsionness is consistent with a growth number conjecture of Mazur (see Theorem \ref{Hida Mazur}). In Section \ref{examples}, we give several examples to illustrate the results of the paper. Finally, in Section \ref{non-com speculation}, we formulate a generalized Conjecture Y for fine Selmer groups over an arbitrary $p$-adic Lie extension which does not necessarily contain the cyclotomic $\Zp$-extension. Here we will show that this conjecture is consistent with the pseudo-nullity conjecture of Coates-Sujatha (see Proposition \ref{torsion psuedo-null non-com}).

\subsection*{Acknowledgement}
The author likes to thank Somnath Jha, Debanjana Kundu, Antonio Lei and Bharathwaj Palvannan for their interest and comments on initial drafts of the paper. He would like to thank the anonymous referee for several helpful comments and suggestions. The author's work is supported by the National Natural Science Foundation of China under Grant No.\ 11771164 and the Fundamental Research Funds for the Central Universities of CCNU under Grant No.\ CCNU20TD002.

\section{Review of some commutative algebra} \label{com alg review}

We collect certain commutative algebraic results that will be required for the discussion of the paper. Throughout this section, $\La$ will always denote a regular local ring. (For our arithmetic purposes, we are usually concerned with regular local rings of the form $\Op\ps{T_1, T_2,..., T_r}$, where $\Op$ is some integral domain which is finite flat over $\Zp$.)
It's a standard fact that $\La$ is therefore a unique factorization domain (cf. \cite[Theorem 20.3]{Matsu}). In particular, it follows from \cite[Theorem 20.1]{Matsu} that every prime ideal of $\La$ of height one is principal.

Recall that a finitely generated $\La$-module $M$ is said to be torsion if for every element $m\in M$, there exists $x\in \La$ such that $x m=0$. Equivalently, this is saying that $\Hom_{\La}(M,\La)=0$. The module $M$ is said to be pseudo-null if the localization $M_\mathfrak{p}$ of $M$ at every prime ideal $\mathfrak{p}$ of height $\leq 1$ is trivial. The latter is equivalent to $\Ext_{\La}^i(M,\La)=0$ for $i=0,1$ (for instance, see \cite[Chap. 5]{Matsu} or \cite[Chap. V, \S1]{NSW}).

We now present a useful lemma (compare with \cite[Section 1.3, Lemme 4]{PR84})

\bl \label{alg lemma}
Let $x$ be an element in $\La$ which is a generator of a prime ideal of $\La$ of height one. Write $\Omega:=\La/x$ for the quotient ring which is also a regular local ring. Then the following statements are valid.
\begin{enumerate}
  \item[$(i)$] If $y$ is another prime element of $\La$ which is coprime to $x$, then $\La/(x,y^n)$ is a torsion $\Omega$-module for every $n\geq 1$.
  \item[$(ii)$] If $M$ is a pseudo-null $\La$-module, then both $M[x]$ and $M/x$ are torsion over $\Omega$.
  \item[$(iii)$] If $M$ is a $\La$-module with $M/x$ being torsion over $\Omega$, then $M$ is torsion over $\La$ and $M[x]$ is torsion over $\Omega$.
\end{enumerate}
\el

\bpf
 We begin proving assertion (i). In fact, we shall establish a slightly stronger assertion: namely, if $z$ is an element of $\La$ which is coprime to $x$, then $\La/(x,z)$ is a torsion $\Omega$-module. Since $z$ is coprime to $x$, it does not lie in the ideal $(x)$. Hence $z+(x)$ is a nonzero element of $\Omega$, and it plainly annihilates $\La/(x,z)$. Therefore, this proves our first assertion.

For the proof of (ii), we shall write $N$ for either $M[x]$ or $M/x$. Since the module $N$ is annihilated by $x$, it may be viewed as a $\Omega$-module. Consider the following spectral sequence
\begin{equation}\label{spec seq}
 \Ext^i_\Omega\big(N, \Ext^j_\La(\Omega, \La)\big)\Rightarrow \Ext^{i+j}_\La(N, \La)
\end{equation}
  (cf. \cite[Exercise 5.6.3]{Wei}). From the $\La$-free resolution
 \[ 0\lra \La\lra \La \lra \Omega \lra 0\]
 of $\Omega$, we see that
 \[\Ext^j_\La(\Omega, \La)=\left\{
                          \begin{array}{ll}
                            \Omega, & \mbox{if $j=1$,} \\
                            0, & \hbox{otherwise.} \\
                                                     \end{array}
                        \right.\]
Therefore, the spectral sequence (\ref{spec seq}) degenerates yielding the isomorphism
\[ \Ext^i_{\Omega}(N, \Omega) \cong \Ext^{i+1}_{\La}(N, \La)\] for $i\geq 0$.
In particular, we have
\[ \Hom_{\Omega}(N, \Omega)=\Ext^{1}_{\La}(N, \La)=0,\]
where the final zero follows from the assumption that $M$ is pseudo-null over $\La$. This in turn implies that $N$ is torsion over $\Omega$ and we have the second assertion.

The final assertion follows from \cite[Section 1.3, Lemme 2]{PR84} or \cite[Corollary 4.13]{LimFine}.
\epf

Now, let $M$ be a finitely generated torsion $\La$-module. By \cite[Proposition 5.1.7]{NSW}, there is a pseudo-isomorphism
\begin{equation}\label{pseudo-iso structure}
\varphi: M \lra \bigoplus_{i\in I}\La/x_i^{n_i},
\end{equation}
where $I$ is a finite indexing set, each $x_i$ is a generator of a prime ideal of height one and $n_i$ is a non-negative integer. Note that the prime ideals $\La x_i$ and integers $n_i$ are determined by the module $M$.

\bl \label{coprime torsion}
Notation as above. Suppose that $x$ is a prime element in $\La$ which is coprime to all the $x_i$'s. Then $M/x$ is a torsion module over the ring $\Omega=\La/x$.
\el

\bpf
Let $P_1=\ker \varphi$, $P_2=\coker\varphi$ and $Q=\mathrm{im}~\varphi$, where $\varphi$ is given as in (\ref{pseudo-iso structure}). From which, we have the following two short exact sequences
\[ 0\lra P_1\lra M\lra Q \lra 0,\]
\[ 0\lra Q\lra \bigoplus_{i\in I}\La/x_i^{n_i} \lra P_2\lra 0.\]
From which, we have
\[  P_1/x\lra M/x\lra Q/x \lra 0,\]
\[  P_2[x] \lra Q/x\lra \bigoplus_{i\in I}\La/(x_i^{n_i},x) \lra P_2/x\lra 0.\]
By Lemma \ref{alg lemma},  the modules $P_1/x$, $P_2[x]$ and $\bigoplus_{i\in I}\La/(x_i^{n_i},x)$ are torsion over $\Omega$. Putting these observations into the above two exact sequences, we see that so is $M/x$.
\epf

\section{Fine Selmer groups of abelian varieties} \label{fineAbvar}

\subsection{Fine Selmer groups}
We begin reviewing the fine Selmer groups of abelian varieties following \cite{CS05a, CS05, WuPhD, Wu, WuTS, Hac, LMu, LimFineDoc}. Fix an odd prime $p$. Let $A$ be an abelian variety defined over a number field $F$. Let $S$ be a finite set of primes of $F$ containing the primes above $p$, the bad reduction primes of $A$ and the infinite primes. Denote by $F_S$ the maximal algebraic extension of $F$ which is unramified outside $S$. For every extension $\mathcal{L}$ of $F$ contained in $F_S$, we write $G_S(\mathcal{L})=\Gal(F_S/\mathcal{L})$, and denote by $S(\mathcal{L})$ the set of primes of $\mathcal{L}$ above $S$.

Let $L$ be a finite extension of $F$ contained in $F_S$. The fine Selmer group of $A$ over $L$ is defined by
\[ R(A/L) =\ker\left(H^1(G_S(L),\Ap)\lra \bigoplus_{v\in S(L)} H^1(L_v, \Ap)\right).\]
We remark that the above definition is independent of the choice of $S$ (see \cite[Lemma 4.1]{LMu}). Just as the classical $p$-primary Selmer group, the fine Selmer group sits in the following analogous short exact sequence
\begin{equation} \label{eqn fine short exact} 0 \lra \mM(A/L) \lra R(A/L) \lra \Zhe(A/L)\lra 0
\end{equation}
where $\mM(A/F)$ is the fine ($p$-)Mordell-Weil group and  $\Zhe(A/F)$  is the fine ($p$-)Tate-Shafarevich group in the sense of Wuthrich \cite{WuTS}. The fine Mordell-Weil group $\mM(A/F)$ is defined to be the subgroup of $A(F)\ot\Qp/\Zp$ consisting of those elements which are mapped to zero in $A(F_v)\ot_{\Zp}\Qp/\Zp$ for all primes $v$ above $p$. It can be shown that $\mM(A/F)$ injects into $R(A/F)$. The fine Tate-Shafarevich group $\Zhe(A/F)$ is then defined to be the cokernel of this injection (see \cite[Section 2]{WuTS} for details).

Let $F_\infty$ be a (not necessarily cyclotomic) $\Zp$-extension of $F$, whose Galois group $\Gal(F_\infty/F)$ will be denoted by $\Ga$. If $\Ga_n$ denotes the unique subgroup of $\Ga$ of index $p^n$, we write $F_n$ for the fixed field of $\Ga_n$. The fine Selmer group of $A$ over $F_\infty$ is defined to be $R(A/F_\infty) = \ilim_n R(A/F_n)$ which comes naturally equipped with a $\Zp\ps{\Ga}$-module structure. The $\Zp\ps{\Ga}$-modules $\mM(A/F_\infty)$ and $\Zhe(A/F_\infty)$ are similarly defined by taking limit of the corresponding objects over the intermediate subfields. We shall write $Y(A/F_\infty)$ for the Pontryagin dual of $R(A/F_\infty)$. Upon taking direct limit of the sequence (\ref{eqn fine short exact}) and following up by taking Pontryagin dual, we obtain
\begin{equation} \label{fine short exact} 0 \lra \Zhe(A/F_\infty)^{\vee}  \lra Y(A/F_\infty) \lra \mM(A/F_\infty)^\vee \lra 0.
\end{equation}
It is not difficult to verify that the modules occurring in the exact sequence are finitely generated over $\Zp\ps{\Ga}$ (for instance, see \cite[Lemma 3.2]{LimFineDoc}). In fact, one  expects more (see \cite{PR00, WuPhD, LimFineDoc}).

\begin{conjecture}[Conjecture Y] Let $A$ be an abelian variety defined over a number field $F$ and $F_\infty$ a $\Zp$-extension of $F$. Then $Y(A/F_\infty)$ is torsion over $\Zp\ps{\Ga}$.
\end{conjecture}

\br
(1) When $F_\infty$ is the cyclotomic $\Zp$-extension, the above conjecture is a consequence of a conjecture of Mazur \cite{Maz} and Schneider \cite{Sch85} on the structure of the classical Selmer groups. The latter is known when $A$ is an elliptic curve over $\Q$ with good reduction at $p$ and $F$ is an abelian extension of $\Q$ (see \cite{K, Kob}).

(2) Suppose that $E$ is an elliptic curve defined over $\Q$ with complex multiplication given by a imaginary quadratic field $K$ at which $p$ split completely in $K/\Q$. Let $K_\infty$ be the $\Zp$-extension of $K$ which is unramified outside one of the primes of $K$ above $p$. Then the validity of Conjecture Y is a consequence of a result of Coates (see \cite[Chap IV., Corollary 1.8]{deS}; also see the recent papers \cite{CKL, Ke}).

(3) Suppose that $E$ is an elliptic curve defined over $\Q$, and $K^{\mathrm{ac}}$ the anti-cyclotomic $\Zp$-extension of an imaginary quadratic field $K$. Then the $\Zp\ps{\Gal(K^{\mathrm{ac}}/K)}$-torsionness of $Y(A/K^{\mathrm{ac}})$ is known in many cases (for instance, see \cite{LV, PW}).
\er

We record two results taken from \cite{LimFineDoc} which serve as support for Conjecture Y, and will play some role in the subsequent discussion of the paper.

\bt \label{control theorem cor}
Let $A$ be an abelian variety defined over a number field $F$. Let $F_{\infty}$ be a $\Zp$-extension of $F$. If $Y(A/F)$ is finite, then $Y(A/F_\infty)$ is torsion over $\Zp\ps{\Ga}$.
\et

\bpf
See \cite[Corollary 3.5]{LimFineDoc}.
\epf

\bp \label{torsion fine TateSha}
 Let $A$ be an abelian variety defined over a number field $F$, and $F_{\infty}$ a $\Zp$-extension of $F$. Suppose that $\Zhe(A/F_n)$ is finite for every $n$. Then  $\Zhe(A/F_{\infty})$ is a cotorsion $\Zp\ps{\Ga}$-module.
\ep

\bpf
See \cite[Proposition 4.1]{LimFineDoc}.
\epf

We now record a corollary of the preceding proposition.

\bc \label{torsion fine TateSha fg}
 Let $A$ be an abelian variety defined over a number field $F$, and $F_{\infty}$ a $\Zp$-extension of $F$. Suppose that the following statements are valid.
 \begin{enumerate}
   \item[$(a)$] $\Zhe(A/F_n)$ is finite for every $n$.
   \item[$(b)$] $A(F_\infty)$ is a finitely generated abelian group.
 \end{enumerate}
 Then $Y(A/F_{\infty})$ is a torsion $\Zp\ps{\Ga}$-module.
\ec

\bpf
From Proposition \ref{torsion fine TateSha}, we see that $\Zhe(A/F_\infty)^\vee$ is torsion over $\Zp\ps{\Ga}$ under hypothesis (a). Hypothesis (b) tells us that $\mM(A/F_\infty)^\vee$ is torsion over $\Zp\ps{\Ga}$. The conclusion follows from these and the short exact sequence (\ref{fine short exact}).
\epf

\subsection{Connection with the pseudo-nullity conjecture of Coates-Sujatha}

We now study the relation between our Conjecture Y and the Conjecture B of Coates-Sujatha \cite[Conjecture B]{CS05}. As a start, we recall their conjecture, which for now is stated for $\Zp^2$-extensions; see Conjecture \ref{Conj B gen} below for the general version.

\begin{conjecture}[Conjecture B]
Let $L_{\infty}$ be a $\Zp^2$-extension of $F$ which contains $F^\cyc$. Then $Y(A/L_\infty)$ is pseudo-null over $\Zp\ps{G}$, where $G=\Gal(L_\infty/F)$.
\end{conjecture}

\br
In \cite{CS05}, they formulated their conjecture under the extra assumption which is their  so-called Conjecture A (see \cite[Conjecture A]{CS05}). In this paper, we do not require this extra hypothesis, and so the above formulation suffices.
\er

Retaining the above notation, we denote by $\Phi(L_\infty/F)$ the set of all $\Zp$-extensions of $F$ contained in $L_\infty$.
For each $\mL\in \Phi(L_\infty/F)$, write $\Ga_{\mL}=\Gal(\mL/F)$ and $H_\mL=\Gal(L_\infty/\mL)$. Fix a topological generator $h_{\mL}$ of $H_\mL$. Then $h_{\mL}-1$ generates a prime ideal of $\Zp\ps{G}$ of height one with
 \[ \Zp\ps{G}/(h_{\mL}-1) \cong \Zp\ps{\Ga_\mL}.\]

We can now establish the following observation as mentioned in the introduction.

\bp \label{torsion psuedo-null Zp2}
 Let $A$ be an abelian variety defined over a number field $F$, and $L_{\infty}$ a $\Zp^2$-extension of $F$ which contains $F^\cyc$. Suppose that Conjecture B is valid for $A$ over $L_\infty$, or in other words, $Y(A/L_\infty)$ is pseudo-null over $\Zp\ps{G}$.
 Then $Y(A/\mL)$ is a torsion $\Zp\ps{\Gal(\mL/F)}$-module for every $\mL\in\Phi(L_\infty/F)$.
\ep

\bpf
From the following commutative diagram
\[   \entrymodifiers={!! <0pt, .8ex>+} \SelectTips{eu}{}\xymatrix{
    0 \ar[r]^{} & R(A/\mL) \ar[d]^{s} \ar[r] &  H^{1}\big(G_{S}(\mL),\Ap\big)
    \ar[d]^{}\\
    0 \ar[r]^{} & R(A/L_\infty)^{H_\mL} \ar[r]^{} & H^{1}\big(G_{S}(L_\infty),\Ap\big)^{H_\mL}
     } \]
we see that $\ker s$ is contained in $H^1(H_\mL,A(L_\infty))$ which is cofinitely generated over $\Zp$. Upon taking Pontryagin dual, we obtain a map
\[ Y(A/F_\infty)_{H_\mL}\lra Y(A/\mL), \]
whose cokernel is finitely generated over $\Zp$. Therefore, for a given $\mL\in\Phi(L_\infty/F)$, whenever $Y(A/L_\infty)_{H_\mL}$ is torsion over $\Zp\ps{\Ga_{\mL}}$, so is $Y(A/\mL)$.

Now, in view of the hypothesis that $Y(A/L_\infty)$ is pseudo-null over $\Zp\ps{G}$, we may apply Lemma \ref{alg lemma}(ii) to conclude that $Y(A/F_\infty)_{H_\mL}$ is torsion over $\Zp\ps{\Ga_{\mL}}$. Combining this with the assertion in the previous paragraph, we have the conclusion.
\epf

We  now state and prove the following.

\bt \label{torsion Zp2}
 Let $A$ be an abelian variety defined over a number field $F$, and $L_{\infty}$ a $\Zp^2$-extension of $F$ which contains $F^\cyc$. Suppose that $Y(A/F^\cyc)$ is torsion over $\Zp\ps{\Gal(F^\cyc/F)}$. Then for all but finitely many $\mL\in\Phi(L_\infty/F)$, Conjecture Y is valid for $Y(A/\mL)$ or, in other words, $Y(A/\mL)$ is torsion over $\Zp\ps{\Ga_\mL}$.
 \et

\bpf
By \cite[Proposition 7.2]{LimFine}, it follows from the $\Zp\ps{\Gal(F^\cyc/F)}$-torsionness of $Y(A/F^\cyc)$ that $Y(A/L_\infty)$ is torsion over $\Zp\ps{G}$. The structure theorem of $\Zp\ps{G}$-modules (cf. \cite[Proposition 5.1.7]{NSW}) then implies that there is a pseudo-isomorphism
\[ Y(A/L_\infty) \sim \bigoplus_{i\in I} \Zp\ps{G}/Q_i^{n_i}\]
of $\Zp\ps{G}$-modules, where $I$ is some finite indexing set and each $Q_i$ is a (principal) prime ideal of $\Zp\ps{G}$ of height one. Since there are only finitely many $Q_i$'s, the element $h_{\mL}-1$ is coprime to these $Q_i$'s for all but finitely many $\mL\in\Phi(L_\infty/F)$. For each of such element $h_{\mL}-1$, it then follows from Lemma \ref{coprime torsion} that $Y(A/L_\infty)_{H_\mL}=Y(A/L_\infty)/(h_\mL-1)$ is torsion over $\Zp\ps{\Ga_{\mL}}$. By a similar argument to that in Proposition \ref{torsion psuedo-null Zp2}, we see that $Y(A/\mL)$ is torsion over $\Zp\ps{\Ga_{\mL}}$ for every such $\mL$. This yields the required conclusion of the theorem.
\epf

We record one case, where we can obtain many cases of validity of Conjecture Y.

\bc \label{torsion Zp2 cor}
 Let $E$ be an elliptic curve defined over $\Q$ and $F$ a finite abelian extension of $\Q$. Let $L_{\infty}$ be a $\Zp^2$-extension of $F$ which contains $F^\cyc$. Then for all but finitely many $\mL\in\Phi(L_\infty/F)$, $Y(E/\mL)$ is cotorsion over $\Zp\ps{\Ga_\mL}$.
 \ec

\bpf
 A well-known theorem of Kato \cite{K} (also see \cite{Kob}) asserts that $R(E/F^\cyc)$ is cotorsion over $\Zp\ps{\Gal(F^\cyc/F)}$. The corollary is now an immediate consequence of this and Theorem \ref{torsion Zp2}.
\epf

We also refer readers to Section \ref{examples} for examples, where $F$ is not necessarily abelian over $\Q$ but Theorem \ref{torsion Zp2} applies.

\section{Fine Selmer groups of elliptic modular forms} \label{finemodform}

As before, $p$ will denote a fixed odd prime.
Let $f$ be a normalized new cuspidal modular eigenform of even weight $k\geq 2$, level $N$ and nebentypus $\epsilon$. Let $\mK_f$ be the number field obtained by adjoining all the Fourier coefficients of $f$ to $\Q$. Throughout, we shall fix a prime $\mathfrak{p}$ of $\mK_f$ above $p$, and let $V_f$ denote the corresponding two-dimensional $\mK_{f,\mathfrak{p}}$-linear Galois representation attached to $f$ in the sense of Deligne \cite{Del}. Writing $\Op=\Op_{\mK_{f,\mathfrak{p}}}$ for the ring of integers of $\mK_{f,\mathfrak{p}}$, we fix a $\Gal(\bar{\Q}/\Q)$-stable $\Op$-lattice $T_f$ in $V_f$. We then set $A_f = V_f/T_f$. Note that $A_f$ is isomorphic to $\mK_{f,\mathfrak{p}}/\Op\oplus \mK_{f,\mathfrak{p}}/\Op$ as $\Op$-modules.

Let $F$ be a finite extension of $\Q$. Denote by $S$ a finite set of primes of $F$ containing those dividing $pN$ and all the infinite primes. Let $F_\infty$ be a $\Zp$-extension of $F$. Following \cite{Jha, JS}, we define the fine Selmer group of $A_f$ over $F_\infty$ to be $\ilim_nR(A_f/F_n)$, where $F_n$ is the intermediate subfield of $F_\infty/F$ with $|F_n:F|=p^n$ and $R(A_f/F_n)$ is defined by
\[ R(A_f/F_n) =\ker\left(H^1(G_S(F_n),A_f)\lra \bigoplus_{v\in S(F_n)} H^1(F_{n,v}, A_f)\right).\]

The Pontryagin dual of $R(A_f/F_\infty)$ is then denoted by $Y(A_f/F_\infty)$.
As before, one can similarly show that $Y(A_f/F_\infty)$ is finitely generated over $\Op\ps{\Ga}$.
The following conjecture is the natural analogue of Conjecture Y for modular forms.

\begin{conjecture}[Conjecture Y$_f$] \label{Conj Yf} Let $A_f$ be defined as above, and $F_\infty$ a $\Zp$-extension of a number field $F$. Denote by $\Ga$ the Galois group $\Gal(F_\infty/F)$. Then $Y(A_f/F_\infty)$ is torsion over $\Op\ps{\Ga}$.
\end{conjecture}

We now present natural analogue of Proposition \ref{torsion psuedo-null Zp2} and Theorem \ref{torsion Zp2} for the fine Selmer group of a modular form.
As a start, we recall the following analogue of Conjecture B which was first studied by Jha \cite{Jha}.

\begin{conjecture}
Let $L_{\infty}$ be a $\Zp^2$-extension of $F$ which contains $F^\cyc$. Then $Y(A_f/L_\infty)$ is pseudo-null over $\Op\ps{G}$, where $G=\Gal(L_\infty/F)$.
\end{conjecture}

As in Section \ref{fineAbvar}, denote by $\Phi(L_\infty/F)$ the set of all $\Zp$-extensions of $F$ contained in $L_\infty$. For each $\mL\in \Phi(L_\infty/F)$, write $\Ga_{\mL}=\Gal(\mL/F)$ and $H_\mL=\Gal(L_\infty/\mL)$.
From which, a similar argument to that in Proposition \ref{torsion psuedo-null Zp2} yields the following.

\bp \label{torsion psuedo-null Zp2 modforms}
  Suppose that $Y(A_f/L_\infty)$ is pseudo-null over $\Op\ps{G}$.
 Then $Y(A_f/\mL)$ is a torsion $\Op\ps{\Gal(\mL/F)}$-module for every $\mL\in\Phi(L_\infty/F)$.
\ep

Similarly, we can establish the following by a similar argument to that in Theorem \ref{torsion Zp2}.

\bt \label{torsion Zp2 modforms}
 Notations as above. Suppose that $R(A_f/F^\cyc)$ is cotorsion over $\Op\ps{\Gal(F^\cyc/F)}$. Then for all but finitely many $\mL\in\Phi(L_\infty/F)$, Conjecture Y$_f$ is valid for $Y(A_f/\mL)$.
 \et

Combining the above with Kato' result \cite{K}, we have the following analogue of Corollary \ref{torsion Zp2 cor}.

\bc
Suppose that $F$ is an abelian extension of $\Q$ and $L_{\infty}$ a $\Zp^2$-extension of $F$ which contains $F^\cyc$. Then for all but finitely many $\mL\in\Phi(L_\infty/F)$, $Y(E/\mL)$ is cotorsion over $\Op\ps{\Ga_\mL}$.
 \ec

We end the section by establishing a control theorem for the fine Selmer groups of elliptic modular forms, which is the analogue to that in \cite[Theorem 3.3]{LimFineDoc} proved for abelian varieties. From now on, we let $S_{fd}$ denote the set of primes in $S$ which do not split completely in $F_\infty/F$. We shall also write $W(L) = W^{\Gal(F_S/L)}$ for any $F \subseteq L \subseteq F_S$. Similarly, for each
$v \in S$, we write $W(L) = W^{\Gal(\bar{F}_v/L)}$ for any $F_v \subseteq L \subseteq \bar{F}_v$.

\bt \label{control theorem}
Let $f$ be a normalized new cuspidal modular eigenform of even weight $k\geq 2$, level $N$ and nebentypus $\epsilon$. Write $A_f$ for the Galois module attached to $f$ defined as above. Let $F_{\infty}$ be a $\Zp$-extension of $F$ with $F_n$ being the intermediate subfield satisfying $|F_n:F|=p^n$. Suppose that either of the following statements is valid.
  \begin{enumerate}
    \item[$(i)$] For every prime $v\in S_{fd}$ and $v_n$ of $F_n$ dividing $v$, then $H^0(F_{n,v_n}, A_f)$ is finite.
    \item[$(ii)$] $K_{f,\mathfrak{p}}\cap \Qp(\mu_{p^\infty})=\Qp$ and $H^0(F_v, A_f)$ is finite for every prime $v\in S_{fd}$.
  \end{enumerate}
  Then the restriction map
\[r_n: R(A_f/F_n) \lra R(A_f/F_{\infty})^{\Ga_n}\]
has finite kernel and cokernel which are bounded independently of $n$.
\et

\br \label{p nmid N remark}
In the case when $p\nmid N$, one automatically has the finiteness of $H^0(F_{n,v_n}, A_f)$ for every $v_n$ above $p$ (see \cite{CSW}).
\er

Before proving Theorem \ref{control theorem}, we establish the following preliminary lemma.

\bl \label{rank leq 1 invariant}
 Let $M$ be a $\Op\ps{\Ga}$-module which is finitely generated over $\Op$. Suppose that either of the following statements is valid.
  \begin{enumerate}
    \item[$(a)$] $M_{\Ga_n}$ is finite for every $n$.
    \item[$(b)$] $K_{f,\mathfrak{p}}\cap \Qp(\mu_{p^\infty})=\Qp$ and $M_{\Ga}$ is finite with $\rank_{\Op}(M)<p-1$.
  \end{enumerate}
  Then the homology group $H_1(\Ga_n, M)$ is finite with order bounded independently of $n$.
\el

\bpf
Suppose that hypothesis (a) is valid. Since $M$ is plainly torsion as a $\Op\ps{\Ga_n}$-module, one therefore has
 \begin{equation} \label{rank} 0 = \rank_{\Op\ps{\Ga_n}}(M) = \rank_{\Op}\big(M_{\Ga_n}\big)- \rank_{\Op}\big(M^{\Ga_n}\big), \end{equation}
where the second equality follows from \cite[Proposition 5.3.20]{NSW}. Combining these observations, we see that each $M^{\Ga_n}$ is finite. This in turn implies that $M^{\Ga_n}$ is contained in $M[p^\infty]$. But since $M$ is finitely generated over $\Op$, the latter is finite, and hence we conclude that $M^{\Ga_n}$ is finite with order bounded independently of $n$. Finally, since $\Ga_n\cong \Zp$, we have the identification $M^{\Ga_n}\cong H_1(\Ga_n,M)$, thus proving the lemma under the validity of the hypothesis (a).

Now suppose that hypothesis (b) is valid. Identify $\Op\ps{\Ga}$ with $\Op\ps{T}$ under a choice of generator of $\Ga$. Since $M_{\Ga}$ is finite, we see that $T$ has to be coprime to the characteristic polynomial of $M$. Since $K_{f,\mathfrak{p}}\cap \Qp(\mu_{p^\infty})=\Qp$, every other cyclotomic polynomial is irreducible over $\Op\ps{T}$. Since such a polynomial has degree $\geq p-1$, it has to be coprime to the characteristic polynomial of $M$. Consequently, $M_{\Ga_n}$ is finite for every $n$. We are therefore in the situation of hypothesis (a), and so the conclusion follows from the above discussion.
\epf

We can now give the proof of Theorem \ref{control theorem}.

\bpf[Proof of Theorem \ref{control theorem}]
Consider the following commutative diagram
\[   \entrymodifiers={!! <0pt, .8ex>+} \SelectTips{eu}{}\xymatrix{
    0 \ar[r]^{} & R(A_f/F_n) \ar[d]^{r_n} \ar[r] &  H^{1}\big(G_{S}(F_n),A_f\big)
    \ar[d]^{h_n} \ar[r] & \displaystyle\bigoplus_{v_n\in S(F_n)} H^1(F_{n,v_n}, A_f) \ar[d]^{g_n=\oplus g_{n, v_n}}\\
    0 \ar[r]^{} & R(A_f/F_\infty)^{\Ga_n} \ar[r]^{} & H^{1}\big(G_{S}(F_\infty),A_f\big)^{\Ga_n}\ar[r] & \displaystyle\left(\bigoplus_{w\in S(F_{\infty})} H^1(F_{\infty,w},A_f)\right)^{\Ga_n}
     } \]
with exact rows. Since $\Ga_n$ has $p$-cohomological dimension 1, the restriction-inflation sequence tells us that $h_n$ is surjective and that $\ker h_n = H^1\big(\Ga_n, A_f(F_\infty)\big)$. It therefore remains to show the finiteness and boundness of $\ker h_n$ and $\ker g_{n}$.

We begin by showing the finiteness and boundness of $\ker g_{n}$. For each $v_n\in S(F_n)$, fix a prime of $F_\infty$ above $v_n$ which is denoted by $w_n$, and write $v$ for the prime of $F$ below $v_n$. Write $\Ga_{w_n}$ for the decomposition group of $w_n$ in $\Ga$. By the Shapiro's lemma and the restriction-inflation sequence, we have
\[ \ker\Big(\bigoplus_{v_n\in S(F_n)} g_{n,v_n} \Big) = \bigoplus_{v_n\in S(F_n)} H^1\Big(\Ga_{w_n}, A_f(F_{\infty,v_n})\Big).\]
 If $v$ is a prime of $F$ below $w_n$ such that $v$ splits completely in $F_\infty/F$, then $\Ga_{w_n}=0$ and so one has $H^1\big(\Ga_{w_n}, A_f(F_{\infty,w_n})\big) =0$. Thus, it remains to consider the primes $v\in S$ which do not split completely in $F_\infty/F$. Since $S$ is a finite set, the number of such possibly nonzero summands $\bigoplus H^1\big(\Ga_{w_n}, A_f(F_{\infty,w_n})\big)$ is therefore finite and bounded independently of $n$.
Hence it remains to show that each $H^1\big(\Ga_{w_n}, A_f(F_{\infty,w_n})\big)$ is finite and bounded independently for those primes lying above $v$ which do not decompose completely in $F_\infty/F$. For this, one just needs to verify that either (a) or (b) of Lemma \ref{rank leq 1 invariant} is valid. We note that hypothesis (a) is a direct consequence of hypothesis (i) of the theorem. It remains to show that hypothesis (ii) of our theorem yields (b) of the said lemma. For this, it suffices to show that $\corank_{\Op}(A_f(F_{\infty,w_n}))<p-1$. We first consider the case when $v$ does not divide $p$. In this setting, $F_{\infty,w_n}$ is the cyclotomic $\Zp$-extension of $F_v$ which is unramified.  Since $v$ divides $N$, $A_f$ cannot be an unramified $\Gal(\bar{F}_v/F_v)$-module and so $A_f$ cannot be realized over $F_{\infty,w_n}$. Therefore, we must have $\corank_{\Op}(A_f(F_{\infty,w_n}))\leq 1$. As the prime $p$ is assumed to be odd, this in turn implies that $\corank_{\Op}(A_f(F_{\infty,w_n}))<p-1$. Now suppose that $v$ divides $p$. It is well-known that $F_v(A_f)$ is a $p$-adic Lie extension of $F_v$ of dimension at least 2 (see \cite{Ri}). It follows from this that $A_f$ cannot be realized over $F_{\infty,w_n}$, and so we have $\corank_{\Op}(A_f(F_{\infty,w_n}))\leq 1$.

We now show that $\ker h_n = H^1(\Ga_n, A_f(F_\infty))$ is finite and bounded independent of $n$. Now since $F_\infty/F$ is a $\Zp$-extension, it must have at least one prime $v$ above $p$ which is ramified in $F_\infty/F$. In view of the discussion in the preceding paragraph, we have $H^0(F_{n,w},A_f)$ being finite, where $w$ is a prime of $F_\infty$ above $v$. This in turn implies that $H^0(G_S(F_n), A_f)$ is finite for every $n$. The desired conclusion is now a consequence of Lemma \ref{rank leq 1 invariant}.
 \epf

Theorem \ref{control theorem} has the following natural corollary.

\bc \label{finite to torsion}
Retain the settings of Theorem \ref{control theorem}. Assume further that $R(A_f/F)$ is finite. Then $Y(A_f/F_\infty)$ is torsion over $\Op\ps{\Ga}$.
\ec

\section{Hida deformations} \label{ordinary sec}

Let us briefly introduce certain notion and facts arising from the work of Hida \cite{Hi86}.
Denote by $\Ga'$ the group of diamond operators for the tower of modular curves $\{Y_1(p^r)\}$. There is a natural identification of $\Ga'$ with $1+p\Zp$ which we denote by $\kappa:\Ga'\stackrel{\sim}{\lra}1+p\Zp$. For an integer $N$ coprime to $p$, we write $h^{\ord}_{\mF}$ for the quotient of the universal ordinary Hecke algebra with conductor $N$ corresponding to an ordinary $\Lambda$-adic eigenform $\mF$. The ring $h^{\ord}_{\mF}$ is a local integral domain which is finite flat over $\Zp\ps{\Ga'}$. In \cite{Hi86}, Hida constructed an irreducible representation
\[ \rho:\Gal(\bar{\Q}/\Q)\lra \mathrm{Aut}_{h^{\ord}_{\mF}}(\mT_{\mF})\]
which is unramified outside $Np$, and where $\mT_{\mF}$ is a finitely generated torsion-free module of generic rank two over $h^{\ord}_{\mF}$. We now impose two standing assumptions on the pair $(\mT_{\mF}, h^{\ord}_{\mF})$.

\begin{enumerate}
  \item[$\mathbf{(H1)}$] The ring $h^{\ord}_{\mF}$ is isomorphic to $\Op\ps{\Ga'}$ for the ring of integers $\Op$ of a finite extension of $\Qp$.
  \item[$\mathbf{(H2)}$] Denote by $\m$ the maximal ideal of $h^{\ord}_{\mF}$. The residual representation $\bar{\rho}\lra \mathrm{Aut}_{h^{\ord}_{\mF}/\m}(\mT_{\mF}/\m \mT_{\mF})$ is absolutely irreducible as a $\Gal(\bar{\Q}/\Q)$-module.
\end{enumerate}
Under $\mathbf{(H2)}$, the module $\mT_{\mF}$ is free over $h^{\ord}_{\mF}$ (see \cite[Section 2, Corollary 6]{MazTil}). Let $\mathfrak{X}_{\mathrm{arith}}(h^{\ord}_{\mF})$ be the set consisting of $\Zp$-algebra homomorphism $\la:h^\ord_{\mF}\lra \bar{\Q}_p$ which satisfies the property that there exists an open subgroup $U$ of $\Ga'$ and a non-negative integer $w$ such that $\la(u) = \kappa^w(u)$ for every $u\in U$. For each of such $\la$, we shall write $w_\la$ for the integer $w$ appearing the above definition. We also write $P_\la$ for the kernel of $\la$ which is a principal prime ideal of $h^{\ord}_{\mF}$ of height one. In fact, identifying $h^{\ord}_{\mF}\cong \Op\ps{T}$ (recall that we are assuming $\mathbf{(H1)}$), the prime ideal $P_\la$ can be viewed as lying above the prime ideal of $\Zp\ps{T}$ generated by $(1+p)^{w_\la}-(1+T)$. We shall write $p_\la$ for a generator of the prime ideal $P_\la$.

By Hida theory, for each $\la\in\mathfrak{X}_{\mathrm{arith}}(h^{\ord}_{\mF})$, there exists a normalized cuspidal eigenform $f_{\la}$ of weight $w_\la+2$ such that $\mT_{\mF}/P_\la\cong T_{f_{\la}}$, where $T_{f_\la}$ is the lattice of the Galois representation attached to $f_{\la}$ as in the sense of Deligne. Here $T_{f_\la}$ is a free $\Op_\la$-module of rank two, where $\Op_{\la}= h^{\ord}_{\mF}/P_{\la}$.

From now on, to simplify notation, we sometimes write $\mT$ for $\mT_{\mF}$ and $\mR$ for $h^{\ord}_{\mF}$. Recall that $\mR\cong\Op\ps{\Ga'}$ by our standing assumption $\mathbf{(H1)}$. Set $\mA = \mT\ot_{\mR}(\mR^{\vee})$. The next lemma is left to the reader as an exercise.

\bl One has $\mA[p_\la] \cong A_{f_\la}$, where $A_{f_\la}$ is the Galois module attached to $f_\la$ as in Section \ref{finemodform}.
\el

Let $F$ be a finite extension of $\Q$. Denote by $S$ a finite set of primes of $F$ containing those dividing $pN$ and all the infinite primes. For a $\Zp$-extension $F_\infty$ of $F$, the fine Selmer group $R(\mA/F_\infty)$ of $\mA$ over $F_\infty$ is defined to be $\ilim_n R(\mA/F_n)$, where $F_n$ is the intermediate subfield of $F_\infty/F$ with $|F_n:F|=p^n$ and $R(\mA/F_n)$ is given by
\[ R(\mA/F_n) =\ker\left(H^1(G_S(F_n),\mA)\lra \bigoplus_{v_n\in S(F_n)} H^1(F_{n,v_n}, \mA)\right).\]

We can now state the following analogue of Conjecture Y for $\mA$.

\begin{conjecture}[Conjecture $\mathcal{Y}$] \label{Conj mY} Retain settings as above. Denote by $Y(\mA/F_\infty)$ the Pontryagin dual of $R(\mA/F_\infty)$. Then $Y(\mA/F_\infty)$ is torsion over $\mR\ps{\Ga}$, where $\Ga=\Gal(F_\infty/F)$.
\end{conjecture}

We now prove a ``hortizontal'' analogue of Theorems \ref{torsion Zp2} and \ref{torsion Zp2 modforms}. In the subsequent discussion, we write $S'_{cd}$ for the set of primes of $S$ which does not divide $p$ and split completely in $F_\infty/F$.

\bt \label{Hida main}
Assume that $\mathbf{(H1)}$ and $\mathbf{(H2)}$ are valid. Suppose that there exists $\eta\in\mathfrak{X}_{\mathrm{arith}}(h^{\ord}_{\mF})$ which satisfies all of the following properties.
\begin{enumerate}
  \item[$(a)$] For every prime $v\in S'_{cd}$, the group $H^0(F_v, A_{f_{\eta}})$ is finite.
  \item[$(b)$] $Y(A_{f_{\eta}}/F_\infty)$ is torsion over $\Op_\eta\ps{\Ga}$.
\end{enumerate}
Then Conjecture $\mathcal{Y}$ is valid, or equivalently, $Y(\mA/F_\infty)$ is torsion over $\mR\ps{\Ga}$. Furthermore, for all but finitely many $\la\in\mathfrak{X}_{\mathrm{arith}}(h^{\ord}_{\mF})$, $Y(A_{f_{\la}}/F_\infty)$ is torsion over $\Op_\la\ps{\Ga}$.
\et

\bpf Consider the following commutative diagram
\[   \entrymodifiers={!! <0pt, .8ex>+} \SelectTips{eu}{}\xymatrix{
    0 \ar[r]^{} & R(A_{f_\eta}/F_\infty) \ar[d]^{r_\eta} \ar[r] &  H^{1}\big(G_{S}(F_\infty),A_{f_\eta}\big)
    \ar[d]^{h_\eta} \ar[r] & \displaystyle\bigoplus_{w\in S(F_\infty)} H^1(F_{\infty,w}, A_{f_\eta}) \ar[d]^{l=\oplus l_{w}}\\
    0 \ar[r]^{} & R(\mA/F_\infty)[p_\eta] \ar[r]^{} & H^{1}\big(G_{S}(F_\infty),\mA\big)[p_\eta]\ar[r] & \displaystyle\left(\bigoplus_{w\in S(F_{\infty})} H^1(F_{\infty,w},\mA)\right)[p_\eta]
     } \]
with exact rows and vertical maps induced by the following short exact sequence
\[ 0 \lra A_{f_\eta}\lra \mA\stackrel{p_\eta}{\lra} \mA\lra 0. \]
We shall first show that the kernel and cokernel of $r_\eta$ are cotorsion over $\Op_\eta\ps{\Ga}$.
To start off, we see that $h_\eta$ is surjective with $\ker h_\eta = \mA(F_\infty)/P_\eta$. Since $\mA(F_\infty)/P_\eta$ is cofinitely generated over $\Op_\eta$, it is cotorsion over $\Op_\eta\ps{\Ga}$. It therefore remains to show that $\ker l$ is cotorsion over $\Op_\eta\ps{\Ga}$. For this, we decompose $l=\oplus_{w\in S(F_\infty)}l_w =\oplus_{v\in S}(\oplus_{w|v}l_w)$ and show that $\ker\big(\oplus_{w|v}l_w\big)$ is cotorsion over $\Op_\eta\ps{\Ga}$ for each $v\in S$. For each $w$, we have $\ker l_w = \mA(F_{\infty,w})/P_\eta$. Now if $v$ is finitely decomposed in $F_\infty$, then the sum $\oplus_{w|v}$ is finite, and so $\ker\big(\oplus_{w|v}l_w\big)$ is cofinitely generated over $\Op_\eta$ for these $v$'s. Now suppose that $v$ splits completely in $F_\infty$. As noted in Remark \ref{p nmid N remark}, $\mA(F_{\infty,w})[p_\eta] = A_{f_\eta}(F_{\infty,w})$ is finite for all $v$ dividing $p$. For those primes not dividing $p$, the finiteness follows from assumption (a). Consequently, $\mA(F_{\infty,w})/P_\eta$ is finite by Lemma \ref{alg lemma}(iii) and is hence annihilated by some powers of $p$. This power of $p$ annihilates $\ker\big(\oplus_{w|v}l_w\big)$. Therefore, we also have $\ker\big(\oplus_{w|v}l_w\big)$ being cotorsion over $\Op_\eta\ps{\Ga}$ for these primes. Hence we have shown that the kernel and cokernel of $r_\eta$ are cotorsion over $\Op_\eta\ps{\Ga}$.

Combining this observation with hypothesis (b) of the theorem, we see that $R(\mA/F_\infty)[p_\eta]$ is cotorsion over $\Op_\eta\ps{\Ga}$. It then follows from Lemma \ref{alg lemma}(iii) that $R(\mA/F_\infty)$ is cotorsion over $\mR\ps{\Ga}$. This proves the first assertion of the theorem.

In view of $\mathbf{(H1)}$, the ring $\mR\ps{\Ga}$ is isomorphic to $\Op\ps{W, T}$, where $1+W$ (resp., $1+T$) corresponds to a topological generator of $\Ga$ (resp., a topological generator of $\Ga'$). By the structure theorem (cf. \cite[Proposition 5.1.7]{NSW}), we then have a pseudo-isomorphism
\[ Y(\mA/F_\infty)^\vee \sim \bigoplus_{i\in I} \mR\ps{\Ga}/Q_i^{n_i}\]
of $\mR\ps{\Ga}$-modules, where $I$ is some finite indexing set and each $Q_i$ is a principal prime ideal of $\mR\ps{\Ga}$ of height one. Since there are only finitely many $Q_i$'s, for all but finitely many $\la\in\mathfrak{X}_{\mathrm{arith}}(h^{\ord}_{\mF})$, we have $Y(\mA/F_\infty)/P_\la$ being torsion over $\Op_{\la}\ps{\Ga}$. From the above discussion, we see that
  \[   Y(\mA/F_\infty)/P_\la \lra Y(A_{f_\la}/F_\infty) \]
  has cokernel which is torsion over $\Op_\la\ps{\Ga}$. It then follows that $Y(A_{f_{\la}}/F_\infty)$ is torsion over $\Op_\la\ps{\Ga}$ for these $\la$'s.
\epf

The finiteness condition $H^0(F_v, A_{f_{\eta}})$ is known to hold in several cases (see \cite[Section 5.2]{HKLR}). In particular, if $\eta$ comes from an elliptic curve, then this is always true, and so we have the following.

\bc \label{elliptic modular hida}
Assume that $\mathbf{(H1)}$ and $\mathbf{(H2)}$ are valid. Suppose that there exists $\eta\in\mathfrak{X}_{\mathrm{arith}}(h^{\ord}_{\mF})$ such that $\mA[P_\eta] \cong \Ep$ for some elliptic curve $E$ with $R(E/F_\infty)$ cotorsion over $\Zp\ps{\Ga}$.

Then $R(\mA/F_\infty)$ is cotorsion over $\mR\ps{\Ga}$. Furthermore, for all but finitely many $\la\in\mathfrak{X}_{\mathrm{arith}}(h^{\ord}_{\mF})$, $R(A_{f_{\la}}/F_\infty)$ is cotorsion over $\Op_\la\ps{\Ga}$.
\ec

\bpf
 It remains to show that $H^0(F_v, \Ep)$ is finite, but this is an immediate consequence of Mattuck's theorem that $E(F_v)$ is finitely generated over $\Op_{F_v}$ (see \cite{Mat}).
\epf

Comparing Theorems \ref{torsion psuedo-null Zp2 modforms} and \ref{Hida main}, one is naturally led to the following question.

\medskip \noindent
\textbf{Question $\mA$.} Is $Y(\mA/F_\infty)$ pseudo-null over $\mR\ps{\Ga}$?
\medskip

To the best knowledge of the author, the structure of $R(\mA/F_\infty)$ does not seem to be well-understood. Even in the context of a cyclotomic $\Zp$-extension, it is still an open question whether $R(\mA/F_\infty)^\vee$ is finitely generated over $\mR$ (see \cite[Conjecture 1]{JS}). There are some recent studies on this structure by Lei-Palvannan \cite{LeiP2} in this direction. However, to the best knowledge of the author, the subject of the fine Selmer group of a Hida deformation over a general $\Zp$-extension does not seem to be covered in the literature.  We can only hope to revisit this subject in a future study.

We end the section by specializing to the case of an imaginary quadratic field, where we explain how a conjecture of Mazur implies the various Conjecture Y's. As a start, we recall the said conjecture of Mazur \cite{Maz83}.

\begin{conjecture}[Growth Number Conjecture] Let $E$ be an elliptic curve defined over $\Q$ and let $K$ denote an imaginary quadratic field. The Mordell-Weil rank of $E$ stays bounded along any $\Zp$-extension of $K$, unless the extension is anticyclotomic
and the root number is negative.
\end{conjecture}

\bt \label{Hida Mazur}
Assume that $\mathbf{(H1)}$ and $\mathbf{(H2)}$ are valid.
 Suppose that there exists $\eta\in\mathfrak{X}_{\mathrm{arith}}(h^{\ord}_{\mF})$ such that $\mA[P_\eta] \cong \Ep$ for some elliptic curve $E$. Let $K$ be an imaginary quadratic field and $K_\infty$ a $\Zp$-extension of $K$. Suppose that all of the following statements are valid.
 \begin{enumerate}
   \item[(a)] If $K_\infty$ is the anticyclotomic $\Zp$-extension, assume further that the root number is positive.
   \item[(b)] The assertion of the growth number conjecture of Mazur is valid for the pair $(E, K_\infty)$. In other words, the Mordell-Weil rank of $E$ stays bounded along the $\Zp$-extension $K_\infty$.
   \item[(c)] The fine Tate-Shafarevich group $\Zhe(E/K_n)$ is finite for every $n$, where $K_n$ is the intermediate subextension of $K_\infty/K$ with $|K_n:K|=p^n$.
 \end{enumerate}
 Then $Y(\mA/K_\infty)$ is torsion over $\mR\ps{\Ga}$. Moreover, for all but finitely many $\la\in\mathfrak{X}_{\mathrm{arith}}(h^{\ord}_{\mF})$, $Y(A_{f_{\la}}/K_\infty)$ is torsion over $\Op_\la\ps{\Ga}$.
\et

\bpf
The growth number conjecture of Mazur implies that $E(K_\infty)$ is a finitely generated abelian group. Therefore, taking assumption (c) into account, we may apply Corollary \ref{torsion fine TateSha fg} to conclude that $Y(E/F_\infty)$ is torsion over $\Zp\ps{\Ga}$. The theorem is now a consequence of a combination of this latter assertion and Corollary \ref{elliptic modular hida}.
\epf

\br
One can of course prove variants of Proposition \ref{torsion psuedo-null Zp2 modforms} and Theorem \ref{torsion Zp2 modforms} for fine Selmer groups of $\mA$. We shall leave the details for the readers to fill in.
\er

\section{Examples} \label{examples}
We give some examples to illustrate our results.

\begin{itemize}
 \item[(i)] Let $E$ be the elliptic curve $y^2 + y = x^3 -x$. Let $p=5$. It is well-known that $\Sel(E/\Q(\mu_{5^{\infty}})) =0$ (cf.\ \cite[Theorem 5.4]{CS10}), where $\Sel(E/\Q(\mu_{5^{\infty}}))$ is the classical Selmer group. Let $F$ be a finite Galois $p$-extension of $\Q(\mu_5)$. Then by \cite[Corollary 3.4]{HM}, $\Sel(E/F^\cyc)$ is cotorsion over $\Zp\ps{\Gal(F^\cyc/F)}$. Since the fine Selmer group $R(E/F^\cyc)$ is contained in $\Sel(E/F^\cyc)$, this in turn implies that $R(E/F^\cyc)$ is cotorsion over $\Z_5\ps{\Gal(F^\cyc/F)}$. Therefore, Theorem \ref{torsion Zp2} applies. In other words, for any given $\Z_5^2$-extension $L_\infty$ of $F$ containing $F^\cyc$, for all but finitely many $\Z_5$-extension $F_\infty$ of $F$ contained in $L_\infty$, we have that $R(E/F_\infty)$ is cotorsion over $\Z_5\ps{\Gal(F_\infty/F)}$.

In view of Conjecture Y, we expect that $R(E/F_\infty)$ is cotorsion over $\Z_5\ps{\Gal(F_\infty/F)}$ for all $\Z_5$-extensions of $F$. But this seems to be out of reach in general. (Note that Proposition \ref{torsion psuedo-null Zp2} cannot apply here, since we do not have a good way of determining pseudo-nullity at our current state of knowledge.)

We however describe how one may obtain the absolute validity of Conjecture Y for the elliptic curve $E$ in question over certain classes of $5$-extensions of $\Q(\mu_5)$. Let $F$ be a finite Galois $5$-extension of $\Q(\mu_5)$ such that every ramified prime $v$ of $F/\Q(\mu_5)$ outside $5$ is neither a split multiplicative reduction prime of $E$ nor a good reduction prime of $E$ with $E(\Q(\mu_5)_v)[5]\neq 0$. By Kida's formula for elliptic curves (cf. \cite[Theorem 3.1]{HM}), $\Sel(E/F^\cyc)$ is finite which in turn implies that $R(E/F^\cyc)$ is finite. It then follows from \cite[Theorem 3.3]{LimFineDoc} that so is $R(E/F)$. Applying Corollary \ref{control theorem cor}, we see that Conjecture Y holds for every $\Z_5$-extension of $F$. Examples of such $F$'s are $\Q(\mu_{5^n}, 5^{-5^m})$.

  \item[(ii)] The above discussion can also be applied to the elliptic curve $E: y^2 + xy = x^3 -x-1$ and $p=7$. Indeed, in this case, one has $\Sel(E/\Q(\mu_{7^{\infty}})) =0$ (cf.\ \cite[Theorem 5.31]{CS10}). For more examples where the discussion in (i) also applies, we refer readers to the tables in \cite{DD} (basically look for those with $\mathcal{L}_E(\sigma)$ being a unit).

\item[(iii)] Even if $\Sel(E/\Q(\mu_{p^\infty}))\neq 0$, there are many numerical examples (for instance, see \cite{DD, PoCom}), where the group $\Sel(E/\Q(\mu_{p^\infty}))$ can be shown to be cofinitely generated over $\Zp$. By virtue of \cite[Corollary 3.4]{HM}, $\Sel(E/F^\cyc)$ is cofinitely generated over $\Zp$ for every finite Galois $p$-extension $F$ of $\Q(\mu_p)$ which in turn implies that $R(E/F)$ is cofinitely generated over $\Zp\ps{\Gal(F^\cyc/F)}$. Hence we can at least apply Theorem \ref{torsion Zp2} for these examples. We mention that the data in \cite{PoCom} also consists of elliptic curves with supersingular reduction at $p$, where the plus-minus Selmer groups in the sense of Kobayashi \cite{Kob} have been verified to be cofinitely generated over $\Zp$. As the fine Selmer group is contained in either of the plus-minus Selmer groups, the fine Selmer group is also cofinitely generated over $\Zp$. This cofinite generation of fine Selmer group is preserved under $p$-base change of fields (for instance, see the proof of \cite[Theorem 5.5]{LMu}). Hence Theorem \ref{torsion Zp2} applies for these examples.

\item[(iv)] Let $E$ be any elliptic curve in $49a$ and $F=\Q(E[7])$. Take $p=7$. We shall show that Conjecture Y is valid for every $\Z_7$-extension of $F$. Indeed, this is plainly true if the $\Z_7$-extension is the cyclotomic $\Z_7$-extension by \cite[Corollary 8]{RT} and \cite[Corollary 3.6]{CS05}. For a non-cyclotomic  $\Z_7$-extension $F_\infty$ of $F$, the compositum of $F_\infty$ and $F^\cyc$ is a $\Z_7^2$-extension of $F$. By \cite[Section 4, Example (a)]{LimPS}, $Y(E/L_\infty)$ is pseudo-null. Therefore, the torsionness of $Y(E/F_\infty)$ follows from this and Proposition \ref{torsion psuedo-null Zp2}. For more of such examples, we refer readers to \cite[Section 4, Example (a)]{LimPS} and \cite[Table 2]{RT}.

\item[(v)] Let $E$ be any elliptic curve in $32a$ and $F=\Q(\sqrt{-43})$. Take $p=3$. It has been verified by Lei-Palvannan that $R(E/L_\infty)$ is pseudo-null (see \cite[Section 8.4]{LeiP}), where $L_\infty$ is the $\Z_3^2$-extension of $F$. Therefore, we may apply Theorem \ref{torsion psuedo-null Zp2} to obtain the validity of Conjecture Y for every $\Z_3$-extension of $F$. \cite[Table 2]{LeiP} provides more examples, where Proposition \ref{torsion psuedo-null Zp2} can be applied to.

\item[(vi)] Let $E$ be the elliptic curve $79a1$ of Cremona's tables given by
\[y^2 + xy +y =x^3 +x^2 -2x.\]
Let $p=3$ and $F=\Q(\mu_3)$.
It has been computed by Wuthrich (see \cite[p. 253] {DD}) that $E(F)$ has $\Z$-rank 1 and the Tate-Sharafevich group $\sha(E/F)[3^\infty]$ is finite. Under these observations, it follows from a similar argument to that in \cite[Theorem 12]{CMc} that one obtains
$H^2(G_S(F),E[3^\infty]) = 0$. By \cite[Lemma 3.2]{Hac}, this in turn implies that $R(E/F)$ is finite. Let $F_\infty$ be any $\Z_3$-extension of $F$. Then Theorem \ref{control theorem cor} applies to yield the torsionness of $R(E/F_\infty)$. As noted in \cite[p. 362]{Jha}, the residual representation $E[3]$ is irreducible. Let $\mF$ be the Hida family associated to $E$. Applying Corollary \ref{elliptic modular hida}, we see that for all but finitely many $\la\in\mathfrak{X}_{\mathrm{arith}}(h^{\ord}_{\mF})$, $R(A_{f_{\la}}/F_\infty)$ is cotorsion over $\Op_\la\ps{\Gal(F_\infty/F)}$.
\end{itemize}

\section{Non-commutative speculation and further remarks} \label{non-com speculation}

In this final section, we formulate an extension of Conjecture Y for a (possibly non-commutative) $p$-adic Lie extension. To simplify the discussion and for better clarity of the ideas behind, we shall restrict our attention to the context of the paper. However, it should be evident that much of the discussion here carries over to broader classes of Galois representations.

We shall let $\bar{A}$ denote either one of the following objects: $A[p^\infty], A_f, \mA$, and $\bar{R}$ the corresponding coefficients rings of these objects; namely, resp., $\Zp, \Op, \mR$. Note that in the context of $\mA$, we still work under the hypotheses $\mathbf{(H1)}$ and $\mathbf{(H2)}$. Let $F_\infty$ be a $p$-adic Lie extension of $F$ such that $G:=\Gal(F_\infty/F)$ is pro-$p$ with no $p$-torsion and that $F_\infty/F$ is unramified outside a finite set of primes of $F$. The fine Selmer group of $\bar{A}$ over $F_\infty$, which is defined similarly as before, now has the structure of a $\bar{R}\ps{G}$-module. We note that the ring $\bar{R}\ps{G}$ is Auslander regular (cf. \cite[Theorem 3.26]{V02}; also see \cite[Theorem A.1]{LimFine}) and has no zero divisors (cf. \cite{Neu}). Therefore, there is a well-defined notion of torsion $\bar{R}\ps{G}$-modules and pseudo-null $\bar{R}\ps{G}$-modules. For our purpose, a $\bar{R}\ps{G}$-module $M$ is said to be torsion (resp., pseudo-null) if $\Ext_{\bar{R}\ps{G}}^i(M,\bar{R}\ps{G})=0$ for $i=0$ (resp., $i=0,1$). The following is then the natural generalization of Conjecture Y for a $p$-adic Lie extension.

\begin{conjecture}[Generalized Conjecture Y] \label{General Conj Y}
For every $p$-adic Lie extension $F_\infty$ of $F$, $Y(\bar{A}/F_\infty)$ is torsion over $\bar{R}\ps{G}$.
\end{conjecture}

For an abelian variety, the above conjecture was formulated in \cite{KL} for a $\Zp^d$-extension of $F$.

\br \label{K-groups remark} We now discuss a context closely related to the generalized Conjecture Y. Had we replaced $\bar{A}$ in the definition of the fine Selmer group by $\Qp/\Zp(-i)$ ($i>0$), where $(-i)$ denotes the $(-i)$th Tate twist, we will obtain the so-called \'etale wild kernel (for instance, see \cite{NQD, NQD02}).
In this context, the author has established the torsionness of this group over \emph{every} $p$-adic Lie extension (see \cite[Proposition 4.1.1]{LimK2n} and \cite[Sections 3.2 and 3.3]{LimEtaleWild} for details). Therefore, this provides some optimism towards the generalized Conjecture Y. \er

A similar argument to that in Theorem \ref{torsion Zp2} yields the following.

\bt \label{torsion Zpd}
Let $d\geq 2$. Suppose that $L_{\infty}$ is a $\Zp^d$-extension of $F$ which contains $F^\cyc$.
Denote by $\Phi_d(L_\infty/F)$ the set of all $\Zp^{d-1}$-extensions of $F$ contained in $L_\infty$.
Assume that $R(\bar{A}/F^\cyc)$ is cotorsion over $\bar{R}\ps{\Gal(F^\cyc/F)}$. Then for all but finitely many $\mL\in\Phi_d(L_\infty/F)$, $R(\bar{A}/\mL)$ is cotorsion over $\bar{R}\ps{\Gal(\mL/F)}$.
 \et

We finally end by providing a conceptual explanation on why the postulation of Conjecture Y over a general $p$-adic Lie extension is plausible. To give this explanation, we now recall the generalized Conjecture B of Coates-Sujatha \cite{CS05} and Jha \cite{Jha}.

\begin{conjecture} \label{Conj B gen}
If $F_\infty$ is a $p$-adic Lie extension of $F$ of dimension $>1$ containing $F^\cyc$, then $Y(\bar{A}/F_\infty)$ is pseudo-null over $\bar{R}\ps{G}$.
\end{conjecture}

\bp \label{torsion psuedo-null non-com}
  Suppose that Conjecture \ref{Conj B gen} is valid for every $p$-adic Lie extension of $F$ of dimension $>1$ containing $F^\cyc$. Then the generalized Conjecture Y is valid for every $p$-adic Lie extension of $F$ of dimension $>1$.
\ep

\bpf
We first make the following remark. Let $G$ be a compact $p$-adic Lie group and $M$ a $\bar{R}\ps{G}$-module. For every open subgroup $G_0$ of $G$, we then have
\[\Ext^i_{\bar{R}\ps{G}}(M,\bar{R}\ps{G}) \cong \Ext^i_{\bar{R}\ps{G_0}}(M,\bar{R}\ps{G_0}). \]
(cf. \cite[Proposition 5.4.17]{NSW}). Therefore, the question of $M$ being torsion (resp., pseudo-null) over $\bar{R}\ps{G}$ is equivalent to $M$ being torsion (resp., pseudo-null) over $\bar{R}\ps{G_0}$.

Now let $F_\infty$ be a $p$-adic Lie extension of $F$ of dimension $>1$. If $F_\infty$ contains $F^\cyc$, then by the assumption of the proposition, $Y(\bar{A}/F_\infty)$ is pseudo-null over $\bar{R}\ps{G}$, and hence is also torsion over $\bar{R}\ps{G}$. Suppose that $F^\cyc$ is not contained in $F_\infty$. Then $F^\cyc\cap F_\infty$ is a finite extension of $F$. In view of the the remark in the first paragraph, we see that neither the hypothesis nor the conclusion is affected if we replace $F$ by $F^\cyc\cap F_\infty$. Hence, upon relabelling, we might as well assume that $F=F^\cyc\cap F_\infty$. Therefore, writing $L_\infty =F^\cyc \cdot F_\infty$, we have $\mathcal{G}:=\Gal(L_\infty/F)\cong Z\times G$, where $Z\cong\Gal(F^\cyc/F)$. Let $z$ be a topological generator of $Z$. Then for every $\bar{R}\ps{\mG}$-module $M$, we have $M_Z = M/(z-1)$ which can be viewed as a $\bar{R}\ps{G}$-module. As $Z$ is central in $\mG$, $M_Z$ may also be viewed as a $\bar{R}\ps{\mG}$-module. In particular,
\[0\lra \bar{R}\ps{\mG}\lra \bar{R}\ps{\mG} \lra \bar{R}\ps{G}\lra 0\]
is a $\bar{R}\ps{\mG}$-free resolution of the $\bar{R}\ps{\mG}$-module $\bar{R}\ps{G}$.
Via a spectral sequence argument similar to that in Lemma \ref{alg lemma}, we obtain
\[  \Ext^i_{\bar{R}\ps{G}}(Y(A/L_\infty)_{Z},\bar{R}\ps{G}) \cong \Ext^{i+1}_{\bar{R}\ps{\mathcal{G}}}(Y(A/L_\infty)_{Z},\bar{R}\ps{\mathcal{G}}). \]
In particular, one has
\[  \Hom_{\bar{R}\ps{G}}(Y(\bar{A}/L_\infty)_{Z},\bar{R}\ps{G}) \cong \Ext^{1}_{\bar{R}\ps{\mathcal{G}}}(Y(\bar{A}/L_\infty)_{Z},\bar{R}\ps{\mathcal{G}}). \]
Since $Y(\bar{A}/L_\infty)$ is pseudo-null over $\bar{R}\ps{\mG}$ by hypothesis, the latter is zero. From which, we see that $Y(\bar{A}/L_\infty)_{Z}$ is torsion over $\bar{R}\ps{G}$. On the other hand, a descent argument as in Theorem \ref{torsion Zp2} yields a map
\[ Y(\bar{A}/L_\infty)_Z\lra Y(\bar{A}/F_\infty), \]
whose cokernel is a quotient of $H^1(Z, \bar{A}(L_\infty))^\vee$, where the latter is plainly finitely generated over $\bar{R}$. It then follows that $Y(\bar{A}/F_\infty)$ is torsion over $\bar{R}\ps{G}$ as required.
\epf

\end{document}